\documentclass{amsart}
\usepackage{amsmath}%
\usepackage{amsfonts}%
\usepackage{amssymb}%
\usepackage{graphicx}

\usepackage{hyperref}
\usepackage{enumerate}
\newtheorem{theorem}{Theorem}
\theoremstyle{plain}

\newtheorem{claim}{Claim}

\newtheorem{conjecture}{Conjecture}
\newtheorem{corollary}{Corollary}

\newtheorem{definition}{Definition}
\newtheorem{example}{Example}

\newtheorem{lemma}{Lemma}

\newtheorem{proposition}{Proposition}
\newtheorem{remark}{Remark}

\numberwithin{equation}{section}

\newtheorem{remarks}{Remarks}

\newcommand{\R}{\mathbb{R}}

\newcommand{\Ric}{\mbox{\normalfont{Ric}}}

\newcommand{\beeq}{\begin{equation}}
\newcommand{\eneq}{\end{equation}}
\newcommand{\beeqs}{\begin{eqnarray*}}
\newcommand{\eneqs}{\end{eqnarray*}}
\newcommand{\besp}{\begin{split}}
\newcommand{\ensp}{\end{split}}
\newcommand{\bepr}{\begin{proof}}
\newcommand{\enpr}{\end{proof}}
\newcommand{\bethr}{\begin{theorem}}
\newcommand{\enthr}{\end{theorem}}
\newcommand{\beths}{\begin{theorem*}}
\newcommand{\enths}{\end{theorem*}}
\newcommand{\becor}{\begin{corollary}}
\newcommand{\encor}{\end{corollary}}
\newcommand{\bere}{\begin{remark}}
\newcommand{\enre}{\end{remark}}
\newcommand{\bers}{\begin{remarks}}
\newcommand{\enrs}{\end{remarks}}
\newcommand{\beres}{\begin{remark*}}
\newcommand{\enres}{\end{remark*}}
\newcommand{\bele}{\begin{lemma}}
\newcommand{\enle}{\end{lemma}}
\newcommand{\beles}{\begin{lemma*}}
\newcommand{\enles}{\end{lemma*}}
\newcommand{\bepro}{\begin{proposition}}
\newcommand{\enpro}{\end{proposition}}
\newcommand{\bepros}{\begin{proposition*}}
\newcommand{\enpros}{\end{proposition*}}
\newcommand{\becl}{\begin{claim}}
\newcommand{\encl}{\end{claim}}
\newcommand{\beex}{\begin{example}}
\newcommand{\enex}{\end{example}}
\newcommand{\beexs}{\begin{example*}}
\newcommand{\enexs}{\end{example*}}
\newcommand{\beco}{\begin{conjecture}}
\newcommand{\enco}{\end{conjecture}}
\newcommand{\becos}{\begin{conjecture*}}
\newcommand{\encos}{\end{conjecture*}}
\newcommand{\bede}{\begin{definition}}
\newcommand{\bedes}{\begin{definition*}}

\newcommand{\hess}{\mbox{\normalfont{Hess}}\,}

\newcommand{\Aand}{\mbox{ and }}
\textheight=\dimexpr
  \ifcase 10pt
    \or 
    \or 
    \or 
    \or 
    \or 
    \or 
    \or 
    \or 
    57\or 
    52\or 
    48\or 
    44\or 
    41\fi 
  \baselineskip+\topskip\relax
\begin{document}
\title[CYLINDRICITY OF COMPLETE EUCLIDEAN SUBMANIFOLDS]{CYLINDRICITY OF COMPLETE EUCLIDEAN SUBMANIFOLDS WITH RELATIVE NULLITY}
\author[Felippe Soares Guimar\~{a}es]{Felippe Soares Guimar\~{a}es}
\address{Felippe Soares Guimar\~{a}es -- Instituto Nacional de Matem\'{a}tica Pura e Aplicada (IMPA) \newline%
\indent Estrada Dona Castorina 110, Rio de Janeiro / Brazil 22640-320}
\email{felippe.guima@gmail.com}%

\author[Guilherme Machado de Freitas]{Guilherme Machado de Freitas}
\address{Guilherme Machado de Freitas -- Politecnico di Torino\newline%
\indent Corso Duca degli Abruzzi, 24, 10129 Turin, Italy}%
\email{guimdf1987@icloud.com}%

\thanks{The first author's research was partially supported by CNPq/Brazil}
\thanks{The second author's research was partially supported by CNPq/Brazil}


\begin{abstract}
We show that a complete Euclidean submanifold with minimal index of relative nullity $\nu_0>0$ and Ricci curvature with a certain controlled decay must be a $\nu_0$-cylinder. This is an extension of the classical Hartman cylindricity theorem.
\end{abstract}

\maketitle

\section{Introduction}\label{intrs}
The simplest examples of isometric immersions $f:M^n\to\R^m$ such that the index of relative nullity is positive everywhere are the $s$-cylinders. The isometric immersion $f$ is said to be an $s$\emph{-cylinder} if there exists a Riemannian manifold $N^{n-s}$ such that $M^n,\ \R^m\Aand f$ have factorizations
\beeqs
M^n=\R^s\times N^{n-s},\ \R^m=\R^s\times\R^{m-s}\Aand f=I\times h,
\eneqs
where $h:N^{n-s}\to\R^{m-s}$ is an isometric immersion and $I:\R^s\to\R^s$ is the identity map. Clearly, in this case the minimal index of relative nullity $\nu_0$ of $f$ is precisely $s$, as long as that of $h$ is zero.

The classical Hartman theorem states that these are the only possible complete examples with nonnegative Ricci curvature.
\bethr[Maltz \cite{MR0643658}]\label{hart}
Let $M^n$ be a complete manifold with nonnegative Ricci curvature and let $f:M^n\to\R^m$ be an isometric immersion with minimal index of relative nullity $\nu_0>0$. Then $f$ is a $\nu_0$-cylinder.
\enthr
The main purpose of this article is to extend the above result to submanifolds with Ricci curvature having a certain controlled decay.
\bethr\label{math}
Let $M^n$ be a complete manifold with
\beeq\label{rice}
\Ric\geq-\left(\hess\psi+\frac{d\psi\otimes d\psi}{n-1}\right)
\eneq
for some function $\psi$ bounded from above on $M^n$ and let $f:M^n\to\R^m$ be an isometric immersion with minimal index of relative nullity $\nu_0>0$. Then $f$ is a $\nu_0$-cylinder.
\enthr
Note that we recover Theorem \ref{hart} from the above by simply taking $\psi$ to be constant.
\bers
(\textrm{i}) In Wylie \cite{wylie}, such a Riemannian manifold satisfying \eqref{rice} was said to be $CD(0,1)$ with respect to the \emph{potential function} $\psi$.\\
(\textrm{ii}) We actually prove a version of Theorem \ref{math} that is more general in two ways. The first is that we can weaken the upper bound on $\psi$ assumption to an integral condition along geodesics, the so-called bounded energy distortion. Secondly the function $\psi$ can be replaced with a vector field $X$. We delay discussing this result until Section \ref{genesec}.
\enrs
\section{Preliminaries}
The main step in the proof of Theorem \ref{math} is Lemma \ref{cyli} below (see Maltz \cite{MR0643658}).
\bele\label{cyli}
Suppose $M^n=\R\times N^{n-1}$ is the Riemannian product of $\R$ and a connected Riemannian manifold $N^{n-1}$, and suppose $f:M^n\to\R^m$ is an isometric immersion mapping a geodesic of the form $\R\times\{q\}$ onto a straight line in $\R^m$. Then $f$ is a 1-cylinder.
\enle
Our result also relies on the fundamental fact that the leaves of the minimum relative nullity distribution of a complete submanifold of $\R^m$ are also complete (cf. Dajczer \cite{MR1075013}).
\bele\label{comp}
Let $M^n$ be a complete Riemannian manifold and let $f:M^n\to\R^m$ be an isometric immersion with $\nu>0$ everywhere. Then, the leaves of the relative nullity distribution are complete on the open subset where $\nu=\nu_0$ is minimal.
\enle
Theorem \ref{hart} follows easily from Lemmas \ref{cyli} and \ref{comp} above together with the Cheeger-Gromoll splitting theorem. Indeed, under the assumptions of Theorem \ref{hart}, Lemma \ref{comp} yields that $M^n$ contains $\nu_0$ linearly independent lines through each point where the index of relative nullity is minimal. By the splitting theorem of Cheeger-Gromoll, $M^n$ is isometric to a Riemannian product $\R^{\nu_0}\times N^{n-\nu_0}$, and Theorem \ref{hart} then follows inductively from Lemma \ref{cyli}.

The proof of our Theorem \ref{math} uses the same ideas above, taking advantage of a recent warped product version of the splitting theorem by Wylie \cite{wylie}. According to this latter result, estimate \eqref{rice} is sufficient to split a complete Riemannian manifold $M^n$ that admits a line into a warped product $\R\times_\rho N^{n-1}$ over $\R$. But since this splitting comes from a line of relative nullity, our goal is to show that the warping function $\rho$ must be constant, and thus $\R\times_\rho N^{n-1}$ is actually a Riemannian product, so that Lemma \ref{cyli} can be applied to conclude the proof. To do this we need to collect geometric information on the behavior of a warped product as above along the line $\R$. For later use, we carry out this study within the broader class of \emph{twisted products} $M^n=\R\times_\rho N^{n-1}$ over $\R$, where $\left(N,h\right)$ is a Riemannian manifold, $\rho:M^n\to\R_+$ the \emph{twisting function}, and $M^n$ is endowed with the metric $g=dr^2+\rho^2h$. If $\rho$ is a function of $r$ only, then we have a \emph{warped product} over $\R$. The following lemma describes how vector fields vary along $\R$.
\bele\label{lccl}
Let $M^n=\R\times_\rho N^{n-1}$ be a twisted product over $\R$. Then
\beeq\label{rr}
\nabla_{\partial_r}\partial_r=0
\eneq
and
\beeq\label{rX}
\nabla_{\partial_r}X=\nabla_X\partial_r=\frac{1}{\rho}\frac{\partial\rho}{\partial r}X
\eneq
for all $X\in\mathfrak{X}(N)$.
\enle
\bepr
Let us write $\rho_r=\rho\left(r,\cdot\right)$ and denote by $N_{\rho_r}$ the Riemannian manifold $N$ endowed with the conformal metric rescaled by $\rho_r^2$. It is straightforward to check that $\nabla$ given by \eqref{rr}, \eqref{rX} and
\beeqs
\nabla_XY=\nabla^{N_{\rho_r}}_XY-\left\langle X,Y\right\rangle\frac{1}{\rho}\frac{\partial\rho}{\partial r}\partial_r
\eneqs
for all $X,\,Y\in\mathfrak{X}(N)$ defines a compatible symmetric connection on $TM$, hence it coincides with the Levi-Civita connection of $M^n$.
\enpr
Next, we use Lemma \ref{lccl} to compute the sectional curvatures along planes containing $\partial_r$.
\bele
Let $M^n=\R\times_\rho N^{n-1}$ be a twisted product over $\R$. Then
\beeq\label{eqhess}
K\left(\partial_r,X\right)=-\frac{1}{\rho}\frac{\partial^2\rho}{\partial r^2}
\eneq
for all unit vector $X\in T_xN$ and all $x\in N^{n-1}$.
\bepr
Differentiating $\left\langle X,X\right\rangle=\rho^2$ twice with respect to $r$ gives
\beeqs
\left\langle\nabla_{\partial_r}\nabla_{\partial_r}X,X\right\rangle+\left\|\nabla_{\partial_r}X\right\|^2=\rho\frac{\partial^2\rho}{\partial r^2}+\left(\frac{\partial\rho}{\partial r}\right)^2.
\eneqs
Using \eqref{rr} and \eqref{rX}, we conclude that
\beeqs
\left\langle R\left(\partial_r,X\right)\partial_r,X\right\rangle=\rho\frac{\partial^2\rho}{\partial r^2},
\eneqs
from which the result follows.
\enpr
\enle
We are now in a position to state and prove our main lemma, in which by a \emph{line of nullity} of a Riemannian manifold $M^n$ we mean a curve $\gamma:\R\to M^n$ such that $\gamma'\left(t\right)\in\Gamma\left(\gamma\left(t\right)\right)$ for all $t\in\R$, where
\beeqs
\Gamma\left(x\right)=\left\{X\in T_xM:R\left(X,Y\right)=0\text{ for all }Y\in T_xM\right\}
\eneqs
is the nullity subspace at $x\in M^n$.
\bele\label{prod}
Let $M^n=\R\times_\rho N^{n-1}$ be a twisted product over $\R$. If $\R\times\{q\}$ is a line of nullity of $M^n$ for some $q\in N^{n-1}$, then $\rho_r=\rho_0$ does not depend on $r$, and hence $M^n$ is actually the Riemannian product $\R\times N_{\rho_0}^{n-1}$.
\enle
\bepr
It follows from \eqref{eqhess} that
\beeqs
\frac{\partial^2\rho}{\partial r^2}\equiv0,
\eneqs
but since the twisting function $\rho$ is positive on the whole real line it must be constant.
\enpr
\section{Proof}\label{proof}
As previously discussed, Lemma \ref{cyli} is at the core of the proof of Theorem \ref{math}, whereas Lemma \ref{prod} is the principle behind its use.
\bepr
We can assume that $\nu_0=1$, since the general case follows easily by induction on $\nu_0$. Take a point $p\in M^n$ where $\nu=1$. It follows from Lemma \ref{comp} that $M^n$ contains a line $l$ through $p$. By the warped product version of the splitting theorem of Cheeger-Gromoll due to Wylie \cite{wylie}, the Riemannian manifold $M^n$ is isometric to a warped product $\R\times_\rho N^{n-1}$ over $\R$, the line $l$ corresponding to $\R\times\{q\}$ for some $q\in N^{n-1}$. Since $l$ is a leaf of the relative nullity foliation, we have in particular that $\R\times\{q\}$ is a line of nullity of $\R\times_\rho N^{n-1}$, and thus, by Lemma \ref{prod}, $\rho_r=\rho_0$ does not depend on $r$ and $\R\times_\rho N^{n-1}$ is actually the Riemannian product $\R\times N_{\rho_0}^{n-1}$.
Hence, we may consider $f:\R\times N_{\rho_0}^{n-1}\to\R^m$, and as $f$ maps $\R\times\{q\}$ onto a straight line in $\R^m$, the result then follows from Lemma \ref{cyli}.
\enpr
\section{Generalization}\label{genesec}
In this section we explain how the result above also has a version for non-gradient
potential fields. Curvature inequality \eqref{rice} has a natural extension to vector fields $X$ and can be regarded as the special case where $X=\nabla\psi$.

Our result in the gradient case assumes boundness of the potential function
$\psi$. While there is no potential function for a non-gradient field, we can still make
sense of bounds by integrating $X$ along geodesics. Let $X$ be a vector field on a
Riemannian manifold $M^n$. Let $\gamma:\left(a,b\right)\to M^n$ be a geodesic that is parametrized
by arc-length. Define
\beeqs
\psi_\gamma\left(t\right)=\int_a^t\left\langle\gamma'\left(s\right),X\left(\gamma\left(s\right)\right)\right\rangle ds,
\eneqs
which is a real valued function on the interval $\left(a,b\right)$ with the property that $\psi_\gamma'\left(t\right)=
\left\langle\gamma'\left(t\right),X\left(\gamma\left(t\right)\right)\right\rangle$. When $X=\nabla\psi$ is a gradient field then $\psi_\gamma\left(t\right)=\psi\left(\gamma\left(t\right)\right)-\psi\left(\gamma\left(a\right)\right)$, in the non-gradient case we think of $\psi_\gamma$ as being the anti-derivative of $X$ along the geodesic $\gamma$. We now recall the notion of `bounded energy distortion', introduced by Wylie \cite{wylie}.
\bede
Let $M^n$ be a non-compact complete Riemannian manifold and $X\in\mathfrak{X}\left(M\right)$ a vector field. Then we say $X$ has \emph{bounded energy
distortion} if, for every point $x\in M^n$,
\beeqs
\limsup_{r\to\infty}\inf_{l\left(\gamma\right)=r}\left\{\int_0^re^{-\frac{2\psi_\gamma\left(\gamma\left(s\right)\right)}{n-1}}ds\right\}=\infty,
\eneqs
where the infimum is taken over all minimizing unit speed geodesics $\gamma$ with $\gamma\left(0\right)=x$.
\end{definition}
In general, $\psi_\gamma$ depends on the parametrization of $\gamma$ only up to an additive
constant, so the notion of bounded energy distortion does not depend on the
parametrization of the geodesic. Also note that if a vector field $X$ has the property
that $\psi_\gamma$ is bounded for all unit speed minimizing geodesics then it has bounded
energy distortion. However, even in the gradient case, bounded energy distortion
is a weaker condition than $\psi$ bounded above.

Our most general cylindricity theorem is the following.
\bethr\label{ngct}
Let $\left(M^n,g\right)$ be a complete manifold with
\beeq\label{rngc}
\Ric\geq-\left(\frac{1}{2}L_Xg+\frac{X^\sharp\otimes X^\sharp}{n-1}\right)
\eneq
for some vector field $X$ with bounded energy distortion and let $f:M^n\to\R^m$ be an isometric immersion with minimal index of relative nullity $\nu_0>0$. Then $f$ is a $\nu_0$-cylinder.
\enthr
In particular, when $X=\nabla\psi$, we conclude that Theorem \ref{math} still holds under the weaker condition that $\psi$ has bounded energy distortion rather than being bounded from above.

By Wylie \cite{wylie}, inequality \eqref{rngc} allows to split $M^n$ as a twisted product $\R\times_\rho N^{n-1}$ over $\R$, provided there is a line. But since Lemma \ref{prod} actually holds for twisted products, the proof of Theorem \ref{ngct} then follows by the same arguments as in Section \ref{proof}.
\bibliographystyle{utphys}
\bibliography{ourbibhart}

\providecommand{\href}[2]{#2}\begingroup\raggedright\begin{thebibliography}{1}

\bibitem{MR0643658}
R.~Maltz, ``Cylindricity of isometric immersions into {E}uclidean space,'' {\em
  Proc. Amer. Math. Soc.} {\bfseries 53} no.~2, (1975) 428--432.

\bibitem{wylie}
W.~Wylie, ``{A warped product version of the Cheeger-Gromoll splitting
  theorem},'' \href{http://arxiv.org/abs/1506.03800}{{\ttfamily
  arXiv:1506.03800 [math.DG]}}.

\bibitem{MR1075013}
M.~Dajczer, {\em Submanifolds and isometric immersions}, vol.~13 of {\em
  Mathematics Lecture Series}.
\newblock Publish or Perish, Inc., Houston, TX, 1990.
\newblock Based on the notes prepared by Mauricio Antonucci, Gilvan Oliveira,
  Paulo Lima-Filho and Rui Tojeiro.

\end{thebibliography}\endgroup
\end{document}